\documentstyle[11pt,amssymb,amsmath,amstext,amsfonts,xypic]{article}
\newcommand{\nc}{\newcommand}

\xyoption{all}   
\CompileMatrices 
\newcommand{\dis}{\displaystyle}

\def\ccirc{{{}_{^{\,^\circ}}}}

\DeclareMathOperator{\Tor}{\mathrm{Tor}}
\DeclareMathOperator{\Ext}{\mathrm{Ext}}

\DeclareMathOperator{\Hom}{\mathrm{Hom}}
\DeclareMathOperator{\RHom}{\mathrm{RHom}}

\newcommand{\limp}{{\underset{^{n}}{\textsl{lim}}}\,
{\textsl{proj}}\;}
\DeclareMathOperator{\REnd}{{\mathrm{REnd}}}
\DeclareMathOperator{\End}{{\mathrm{End}}}
\DeclareMathOperator{\Spec}{\mathrm{Spec}}
\DeclareMathOperator{\pr}{pr}
\newcommand{\opp}{{\operatorname{op}}}
\newcommand{\op}{\operatorname}
\nc{\can}{{\mathsf{can}}}

\renewcommand{\int}{_{_{\sf{internal}}}}
\renewcommand{\o}{\otimes}

\newcommand{\lo}{{\,\overset{L}{\otimes}\,}}

\newcommand{\DGMod}{\op{DGM}}

\newcommand{\r}{{\mathsf{R}}}

\newcommand{\cone}{{\sf{Cone}}}
\newcommand{\rdeform}{{\mathfrak{Deform}}}
\newcommand{\ra}{{\r\aa}}
\newcommand{\fm}{{\mathfrak{m}}}

\newcommand{\rop}{{\r\aa}\otimes{\r\aa}^{\opp}}
\newcommand{\ror}{{\r\o_\oo\r}}

\nc{\der}{{\mathsf Der}} 
\nc{\Def}{{\mathfrak{Deform}}}
\renewcommand{\AA}{{A}}
\renewcommand{\aa}{{\mathsf a}} 

\nc{\tilb}{\til{\mathsf b}}

\nc{\fU}{{\mathfrak{U}}}
\nc{\Ind}{{\mathtt Ind}}

\nc{\Coind}{{\mathtt Coind}} 
\nc{\Ker}{{\mathtt Ker}}
\nc{\im}{{\mathtt Im}} 
\nc{\Coker}{{\mathtt Coker}}
\nc{\dirlim}{\underset{\rightarrow}{{\mathtt lim}}}
\nc{\invlim}{\underset{\leftarrow}{{\mathtt lim}}} 
\nc{\CN}{{\mathcal N}}
\nc{\La}{{\boldsymbol{\Lambda}}}
\nc{\Th}{\Theta}
\nc{\mto}{\longmapsto}
\nc{\lth}{{l} t} 
\nc{\ao}{{A}}

\nc{\qis}{\stackrel{\tt{qis}}{\,\sim\,}}
\nc{\qisto}{\stackrel{\tt{qis}}{\too}}

\nc{\zz}{{\mathfrak{z}}} 
\nc{\til}{\widetilde}
\nc{\fI}{{\mathfrak{I}}}

\nc{\X}{{\mathbb X}} \nc{\Y}{{\mathbb Y}} 
\nc{\YY}{{\fr\Y}}

\nc{\et}{{\boldsymbol{\eta}}}
\nc{\pp}{{\boldsymbol{\pi}}}
\nc{\ten}{{\otimes}} 
\nc{\tenl}{\overset{{\mathtt L}}\ten}
\nc{\map}{\longrightarrow}

\nc{\triv}{\mathsf{triv}} 
\nc{\bs}{\bigskip\\} 
\nc{\ms}{\smallskip\\}
\newtheorem{theorem}[equation]{Theorem}
\newtheorem{prop}[equation]{Proposition}

\newtheorem{lemma}[equation]{Lemma}

\makeatletter
\renewcommand{\subsection}{\@startsection{subsection}{2}{0pt}{-3ex
plus -1ex minus -0.2ex}{-2mm plus -0pt minus
-2pt}{\normalfont\bfseries}} \makeatother

\numberwithin{equation}{subsection}

\def\proof{\smallskip\noindent {\it Proof.  \ }}

 \def\beq{\begin{equation}\label}

 \def\eeq{\end{equation}}

\newcommand{\iso}{{\,\stackrel{_\sim}{\to}\,}}
\newcommand{\cd}{\!\cdot\!}

\def\vi{{\sf{(i)}\,\,}}
\def\vii{{\sf{(ii)}\,\,}} 
\def\viii{{\sf{(iii)}\,\,}}

\newcommand{\tooo}{{\;{-\!\!\!-\!\!\!-\!\!\!-\!\!\!\longrightarrow}\;}}
\nc{\vp}{\varpi} 
\nc{\supp}{{\mathtt{supp}^{\,}}}
 \nc{\starb}{\star}
\nc{\starg}{\star} 
\nc{\stara}{\star_{\!\!_{{\mathtt{act}}}}}

\newcommand{\id}{{\operatorname{Id}}}

\newcommand{\m}{{\mathfrak{m}}}

\newcommand{\into}{\hookrightarrow}
\newcommand{\onto}{\,\twoheadrightarrow\,}

\newcommand{\too}{\;\longrightarrow\;}

\DeclareMathOperator{\sym}{\mathrm{Sym}}

\newcommand{\taul}{{{{}^{\tau}}\!<}}
\newcommand{\taug}{{{}^\tau\geq}}

 \newcommand{\bimod}{\mbox{-}{\tt{bimod}}}

 \newcommand{\oo}{{\mathcal{O}}}

\newcommand{\bm}{{\mathbf{m}}}
\newcommand{\dg}{${\mathtt{dg}}$}

\newcommand{\odg}{${\mathcal{O}}$-${\mathtt{dg}}$}

\newcommand{\qed}{$\quad\square$}

 \def\Z{{\mathbb{Z}}} 
 \def\I{{\Omega}}

\def\k{{\Bbbk}}

\newcommand{\ta}{(\aa\lo_{_{\AA}}\aa)^{\taug -1}}

\def\dim{{\mathsl{dim}\,}}

\def\zz{{\mathcal{Z}}}

\newcommand{\ab}{{\hskip 8mm}}

 \def\triv{{\mathsf {triv}}}

\newcommand{\sset}{\subset}

\newcommand{\pa}{\partial}

\def\square{\hbox{\vrule\vbox{\hrule\phantom{o}\hrule}\vrule}}

\begin{document}

\centerline{{\Large{\bf On Deformations of Associative Algebras}}}
 \bigskip

\centerline{{\sc Roman Bezrukavnikov and Victor 
Ginzburg}\footnote{Both authors are partially supported by the NSF.}}
\vskip 6mm

\begin{abstract}{\footnotesize
\noindent In a classic paper, Gerstenhaber showed that first
order deformations of an associative $\k$-algebra $\aa$ are controlled
by the second Hochschild cohomology group of $\aa$. 
More generally,
any   $n$-parameter  first
order deformation of $\aa$ gives,
due to commutativity of the cup-product on  Hochschild cohomology,
a graded algebra morphism $\sym^\bullet(\k^n) \to \Ext^{2\bullet}_{\aa\bimod}(\aa,\aa)$.
We prove that any
extension of the
 $n$-parameter  first
order deformation of $\aa$ to an {\em infinite order} formal
deformation
provides a canonical `lift' of 
the graded algebra morphism above
 to a  \dg-algebra morphism $\sym^\bullet(\k^n) \to 
\RHom^{\bullet}_{}(\aa,\aa),$
where the Symmetric algebra $\sym^\bullet(\k^n)$ is viewed as a
\dg-algebra 
(generated by the vector space $\k^n$ placed in degree 2)
equipped with zero differential.
}
\end{abstract}

%sss

\section{Main result}
\setcounter{equation}{0}
\subsection{}\label{setting} Let $\k$ be a
field of characteristic zero and 
write $\o=\o_\k,\,\Hom=\Hom_\k,$ etc. 
Given a  $\k$-vector space $V$, let $V^*=\Hom(V,\k)$ denote the
dual space.

We will work with unital associative $\k$-algebras,
to be referred as `algebras'. Given such an algebra $B$,
we write $\bm_B: B\o B\to B$
for the corresponding multiplication
map, and put  $\I_B:=\Ker(\bm_{B})\sset B\o B$. This
is a  $B$-bimodule which is
free as a  {\em right} $B$-module;  in effect,
$\I_{B}\simeq (B/\k)\o B$ is a free
right $B$-module generated by the
subspace $B/\k\sset \I_{B}$ formed by the elements
$b\o 1-1\o b,\,b\in B$.

Fix a finite dimensional vector space $T,$ and let
$\oo=\k\oplus T^*$ be the commutative local $\k$-algebra
with unit $1\in\k$ and with maximal ideal $T^*\sset \oo$ such that
 $(T^*)^2=0.$ Thus, $\oo/T^*=\k$.
The algebra $\oo$ is Koszul and one has
a canonical isomorphism
$\Tor^\oo_1(\k,\k)\cong T^*.$

We are interested in multi-parameter (first order)
 deformations of a given algebra
  $\aa$. Specifically, 
by an  $\oo$-deformation of $\aa$
we mean a
 free $\oo$-algebra  $A$ (that is, $\oo$ is a {\em central}
subalgebra in $A,$ and $A$ is a 
 free $\oo$-module) equipped with
a $\k$-algebra isomorphism $\psi: A/ T^*\cd A\cong \aa$.
Two   $\oo$-deformations, $(A,\psi)$ and $(A',\psi'),$ 
 are said to be {\it
equivalent} if there is an $\oo$-algebra isomorphism
$\varphi: A\iso A'$ such that its reduction modulo the maximal ideal
induces the identity map $\id_\aa: \aa\stackrel\psi\simeq
A/ T^*\cd A\stackrel\varphi\too A'/ T^*\cd A'\stackrel{\psi'}\simeq\aa$.

Let  $(A,\psi)$ be an $\oo$-deformation of $\aa$.
Reducing  each term of the short exact sequence
$0\to \I_{A}\to A\o A\to A\to 0$
(of free right $A$-modules) modulo $ T^*$
on the right, one obtains the following
 short exact sequence of left $A$-modules
\beq{classI}
0\to 
\I_{A}\o_\oo
\k\to A\o\aa\to\aa\to0.
\eeq

Next, we reduce  modulo $ T^*$
on the left, that is, apply the functors $\Tor^\oo_\bullet(\k,-)$
with respect to the left $\oo$-action.
We have
 $\Tor^\oo_1(\k,A\o\aa)=0.$
Further, since multiplication by $ T^*$ annihilates $\aa$,
we get
$\Tor^\oo_1(\k,\aa)=\aa\o\Tor^\oo_1(\k,\k)=\aa\o T^*.$
Thus, the end of the long  exact sequence of Tor-groups
corresponding to the short exact sequence
\eqref{classI} reads
\beq{I_A}
0\map\aa\o T^*\stackrel{\nu}\map\k\otimes_{_\oo}\I_{A}\otimes_{_\oo}\k
\stackrel{u}\map\aa\ten\aa
\stackrel{\bm_\aa}\map\aa\map 0.
\eeq
This is an exact sequence of  $\aa$-bimodules;
the map
 $\nu: \aa\o T^*=\Tor^\oo_1(\k,\aa)\to
\k\o_\oo(\I_{A}\otimes_{_\oo}\k)$  in \eqref{I_A}  is the boundary map which
 is easily seen to
be induced by the assignment
$a\o t\mapsto ta\o1-1\o at\in \I_{A}$, for any $a\in 
A$ and $t\in T^*.$ The map $u$ 
 is induced by the imbedding $\I_{A}\into A\o A$.

\smallskip
\noindent
{\sc{Interpretation via  noncommutative
geometry.}}\enspace 
 For any associative algebra
$A$, the bimodule $\I_A$ is called the bimodule 
of {\em noncommutative 1-forms} for $A$, and there is a
geometric
interpretation of \eqref{I_A} as follows.

Let $J\sset A$ be
any two-sided ideal, and put $\aa:=A/J$.
There is a canonical short exact sequence
of $\aa$-bimodules (cf. \cite[Corollary 2.11]{CQ}),
\beq{CQ}
0\too J/J^2 \stackrel{d}\too
\aa\otimes_A\I_A\otimes_A\aa\too \I_\aa\too 0.
\eeq
Here, the  map 
$J/J^2\to\aa\otimes_A\I_A\otimes_A\aa$ is induced by  restriction to $J$ of
the de Rham differential $d: A\to \I_ A$, cf. \cite{CQ}.
The above  exact sequence  may be thought
of as a noncommutative analogue of
the {\em conormal exact sequence} of a subvariety.

We may splice \eqref{CQ} with the tautological extension \eqref{classI},
the latter
tensored by $\aa$ on both sides.
 Thus, we obtain
the following exact sequence of $\aa$-bimodules
\beq{diag}
0\to J/J^2\stackrel{d}\too
\aa\o_A\I_ A\otimes_A\aa\too\aa\otimes\aa\stackrel{{\mathbf{m}}_\aa}\too
\aa\to 0.
\eeq

Let $\Ext_{\aa\bimod}^i(-,-)$ denote the $i$-th
Ext-group in  $\aa\bimod$, the abelian
 category of $\aa$-bimodules.
 The group
$\Ext_{\aa\bimod}^2(\aa, J/J^2)$ classifies 
 $\aa$-bimodule extensions of $\aa$ by $J/J^2$.
The class of  the extension in \eqref{diag} may be viewed
as a noncommutative version of  Kodaira-Spencer class.

We return now to the special case where $A$ is an $\oo$-deformation of
an algebra $\aa$. 
In this case, we have $\aa=A/J$ where $J=\aa\otimes T^*$ and, moreover, $J^2=0$. Thus,
$J/J^2=\aa\otimes T^*$, and the long exact sequence in
\eqref{diag} reduces to
\eqref{I_A}. Let
$${\sf{deform}}(A,\psi)\in\Ext_{\aa\bimod}^2(\aa,\aa\o T^*)=
\Hom(T,\,\Ext_{\aa\bimod}^2(\aa,\aa)).
$$
be the  class  of  the corresponding extension.

The following theorem is an invariant and multiparameter generalization
of the classic result due to Gerstenhaber \cite{G2}.

\begin{theorem}\label{Ger_def} 
The map assigning
the class ${\sf{deform}}(A,\psi)\in\Hom(T,\,\Ext_{\aa\bimod}^2(\aa,\aa))$
 to
an $\oo$-deformation  $(A,\psi)$
provides
a canonical bijection between the set of equivalence
classes of  $\oo$-deformations of 
$\aa$
and the vector space $\Hom(T,\,\Ext^2_{\aa\bimod}(\aa,\aa))$.\qed
\end{theorem}

Gerstenhaber worked in more
down-to-earth terms involving  explicit cocycles.
To make a link with
Gerstenhaber's  formulation, observe that, for any deformation
$(A,\psi),$
the composite $A\onto A/T^*\cd A\stackrel{\psi}\iso\aa$
 may be lifted (since
$A$ is free over $\oo$) to
an $\oo$-module isomorphism
$
A\cong\aa\o\oo
=\aa\o(\k\oplus T^*)=
\aa\bigoplus (\aa\o T^*)$ that
 reduces to $\psi$ modulo $ T^*.$
Transporting the multiplication map
on $A$ via this isomorphism,  we see that giving a   deformation
amounts
to giving an associative truncated `star product'
\beq{star1}
a\star a'= a\cdot a' + \beta(a,a'),\quad
\beta\in \Hom(\aa\otimes \aa, \aa\o T^*)=
\Hom\bigl(T,\,\Hom(\aa\otimes \aa,\aa)\bigr).
\eeq
The associativity condition gives a constraint on $\beta$ 
 saying that  $\beta$ is a  Hochschild 2-cocycle 
(such that $\beta(1,x)=0$).
Changing the  isomorphism
$A\simeq \oo\otimes \aa$ has the effect of replacing  $\beta$ 
by a cocycle in the same cohomology class.

One can show that $\beta={\sf{deform}}(A,\psi),$ i.e.,
the class of the 2-cocycle $\beta$ in the Hochschild cohomology
group
$\Ext_{\aa\bimod}^2(\aa,\aa\o T^*)$
represents the class of the extension in \eqref{I_A}.
Thus, our cocycle-free construction is  equivalent to the one
given by Gerstenhaber.

\subsection{}\label{inf} Next, we consider the total 
Ext-group
$\Ext^\bullet_{\aa\bimod}(\aa,\aa)=\bigoplus_{i\geq 0}
\Ext^i_{\aa\bimod}(\aa,\aa).$
This is a graded  vector space that comes equipped with
 an  associative algebra
structure given by   Yoneda product.
Another fundamental result due to
Gerstenhaber \cite{G1}
says
\begin{theorem}\label{Ger_com}
 The Yoneda product on $\Ext^\bullet_{\aa\bimod}(\aa,\aa)$
is (graded) commutative.\qed
\end{theorem}

 In view of this result,
 any linear map
$T\map \Ext_{\aa\bimod}^2(\aa,\aa),$
of vector spaces,
 can be uniquely extended,
due to  commutativity of the algebra  
$\Ext_{\aa\bimod}^{2\bullet}(\aa,\aa)$,
to a graded algebra homomorphism
$\sym(T[-2])\to \Ext_{\aa\bimod}^{2\bullet}(\aa,\aa),$
 where $\sym(T[-2])$ denotes the commutative graded algebra
freely generated by the vector space $T$
placed in degree 2.
We conclude that any   $\oo$-deformation of
$\aa$
gives rise, by Theorem \ref{Ger_def},
to  a graded algebra homomorphism
\beq{rise}
{\mathsf {deform}}: \
\sym(T[-2])\to \Ext_{\aa\bimod}^{2\bullet}(\aa,\aa).
\eeq

The present paper is concerned with the problem of `lifting'  this
morphism 
to the level of derived categories. 
Specifically,   we
  consider
the \dg-algebra $\RHom_{\aa\bimod}(\aa,\aa)$, see
Sect. \ref{DG} below, and also
view the graded algebra $\sym(T[-2])$ as a \dg-algebra with
trivial differential.  We are interested in
lifting the graded algebra map \eqref{rise}
to a \dg-algebra map $\sym(T[-2])\to\RHom_{\aa\bimod}(\aa,\aa)$.

To this end, one has
to consider {\em infinite order} formal deformations of $\aa$.
Thus, we now let $\oo$ be a {\em formally smooth}
local $\k$-algebra with maximal ideal $\m$ such that $\oo/\m=\k$.
We assume $\oo$ to be complete in
the $\m$-adic topology, that is,
 $\oo\cong\limp\oo/\m^n$.
 The (finite dimensional)
$\k$-vector space  $T:=(\m/\m^2)^*$ may be viewed
as the tangent space
to $\Spec\oo$ at the base point and
one has a {\em  canonical} isomorphism $\oo/\m^2=\k\oplus T^*.$
The algebra $\oo$ is {\em noncanonically}  isomorphic to
$\k[[T]],$ the algebra of formal power series
on the vector space~$T.$

Let  $A$ be a complete topological
$\oo$-algebra, $A\cong\limp A/\m^n A$,
such that, for any $n=1,2,...,$
the quotient $A/\m^n A$ is a free
$\oo/\m^n$-module. Given  an algebra $\aa$
and an
algebra isomorphism
$\psi: \aa\iso A/\m A$, we say that the pair $(A,\psi)$ is an infinite order
formal $\oo$-deformation of $\aa$.

Clearly, reducing an infinite order deformation  modulo
$\m^2$, one obtains a  first order  $\oo/\m^2$-deformation of 
$\aa$.
The main result of this paper reads

\begin{theorem}[Deformation formality]
\label{General_formality1} Any infinite order
formal  $\oo$-deformation $(A,\psi)$ of an associative
algebra $\aa$ 
provides a canonical lift
of the graded algebra morphism
\eqref{rise}, associated
with the corresponding first order $\oo/\m^2$-deformation,
to a \dg-algebra morphism 
$\rdeform: \sym(T[-2])\to \RHom_{\aa\bimod}(\aa,\aa),$
see \S\ref{DG} for explanation.
\end{theorem}

Observe   that  Theorem \ref{General_formality1}
says, in particular, that one  can map
a basis of the vector space
${\mathsf {deform}}(T[-2])\sset \Ext_{\aa\bimod}^2(\aa,\aa)$
to a set of {\em pairwise commuting} elements
in $\RHom_{\aa\bimod}(\aa,\aa).$ Thus, the above theorem
may be seen as a (partial) refinement of 
Gerstenhaber's Theorem \ref{Ger_com}. Yet, our approach to
 Theorem \ref{General_formality1}
is {\em totally different} from Gerstenhaber's
proof of his theorem; indeed, we are unaware of
any connection between the commutativity
resulting from  Theorem \ref{General_formality1}
and  the Gerstenhaber
{\em brace operation} on Hochschild cochains that plays a crucial
role in the proof of Theorem~\ref{Ger_com}. 
This `paradox' may be resolved, perhaps, by observing
that the notation $\RHom_{\aa\bimod}(\aa,\aa)$
stands for a {\em quasi-isomorphism class}
of DG algebras, see \S\ref{DG} below. Yet, the very notion of
 commutativity of elements of
$\RHom_{\aa\bimod}(\aa,\aa)$ only makes sense  after one
picks a concrete  DG algebra in that
quasi-isomorphism class.  Thus, the commutativity statement
resulting from Theorem \ref{General_formality1} implicitly
involves a particular DG algebra model
for $\RHom_{\aa\bimod}(\aa,\aa)$.
Now, the point is 
that  the model that we are using as well as
our construction of the morphism $\rdeform$ will both involve the
full  {\em infinite order} deformation $(A,\psi)$, 
i.e., the full $\oo$-algebra structure on $A$, 
and not only
the `first order' deformation $A/\m^2 A$.
On the contrary, the statement of Gerstenhaber's  Theorem \ref{Ger_com}
is independendent of  the choice of  a  DG algebra model;
also, the construction of the map ${\mathsf{deform}}$ in \eqref{rise}
involves the first order deformation $A/\m^2 A$
only.

\smallskip
\noindent{\em Remark.}$\enspace$
 Theorem  \ref{General_formality1} was applied in
\cite{ABG} to certain natural deformations
of quantum groups at a root of unity.

\smallskip
\noindent{\em Acknowledgements.}$\enspace$ We would like to thank
Vladimir Drinfeld for many interesting discussions which motivated,
in part, a key construction of this paper.

\section{Generalities}
\subsection{Reminder on dg-algebras and dg-modules.}\label{DG}
Given an integer $n$ and a graded vector space $V=\bigoplus_{i\in\Z}\,V^i,$
we write $V_{<n}:=\bigoplus_{i<n}\,V^i$.
Let $\,[n]\,$ denote the shift functor in the derived
category, and also the grading shift by $n$, i.e., $(V[n])^i:= V^{i+n}.$

Let $B=\bigoplus_{i\in\Z}\,B^i$ be a dg-algebra.
We write $\DGMod(B)$ for the homotopy category of all
left dg-modules
$M=\bigoplus_{i\in\Z}\, M^i$ over $B$ (with
differential $d: M^\bullet\to M^{\bullet+1}$),
 and $D(B):=D(\DGMod(B))$ for the corresponding derived category
obtained by localizing at quasi-isomorphisms.
A $B$-bimodule is the same thing as a
left module over $B\otimes B^\opp,$ where
$B^\opp$ stands for the opposite algebra.
Thus, we write
$D(B\otimes B^\opp)$ for
the derived category of \dg-{\em bimodules} over $B$.

Given two objects  $M,N \in D(B),$ for any integer $i$
we put $\Ext_{_{B}}^i(M,N):=
\Hom_{D(B)}(M,N[i]).$
The graded space
$\Ext_{_{B}}^\bullet(M,M)$
$=\bigoplus_{j\geq
0}\,\Ext_{_{B}}^j(M,M)$
 has a natural algebra structure, via composition.

Given   an exact triangle
$\Delta: K\to M\to N$, in $D(B),$
we write 
$\pa_\Delta: N\to K[1]$ for the corresponding boundary
morphism. Thus, 
$\pa_\Delta\in\Hom_{D(B)}(N,K[1])=\Ext^1_B(N,K).$

For a \dg-algebra $B=\bigoplus_{i\leq 0}\,B^i$ 
concentrated in {\em nonpositive} degrees, the triangulated
category $D(B)$ has a standard $t$-{\em structure}
$\dis(D^{\taul 0}(B), D^{\taug 0}(B))$
where $D^{\taul 0}(B)$, resp. $D^{\taug 0}(B)$, is a full subcategory
of $D(B)$ formed by the objects with vanishing cohomology
in degrees $\geq 0$, resp., in degrees $< 0$, cf. \cite{BBD}.
Write $\dis D(B)\to D^{\taul 0}(B),\, M\mapsto M^{\taul 0},$
resp., $\dis D(B)\to D^{\taug 0}(B),\, M\mapsto M^{\taug 0},$
for the corresponding truncation functors.
Thus, for any object $M\in D(B)$, there is a canonical
exact triangle $\dis M^{\taul 0}\to M\to M^{\taug 0}.$
A triangulated
functor $F: D(B_1)\to D(B_2)$ between two such categories
is called $t$-{\em exact} if
it takes $ D^{\taul 0}(B_1)$ to $D^{\taul 0}(B_2),$
 and $ D^{\taug 0}(B_1)$ to $D^{\taug 0}(B_2)$.

 An object $M \in \DGMod(B)$ is said to be
{\it projective} if it belongs to  the smallest full
subcategory of $\DGMod(B)$ that contains the rank one
\dg-module $B$,   and which
is closed under taking
mapping-cones and infinite direct sums. 
Any object of
$\DGMod(B)$ is quasi-isomorphic to a projective
object, see [Ke] for a proof. (Instead of projective objects, one can
use {\em semi-free objects} considered e.g.
in \cite[Appendices A,B]{Dr}.)

Given $M\in \DGMod(B)$, 
choose a quasi-isomorphic projective object $P\in \DGMod(B)$
and write $\Hom'(P, P[n])$ for the space of $B$-module maps
$P\to P$ which shift the grading by $n$ (but do not necessarily commute
with the differential $d$).
The graded vector space
$\bigoplus_{n\in\Z}\,\Hom'(P, P[n])$
has a natural  algebra structure given by composition.
Super-commutator with the differential $d\in \Hom_\k(P, P[-1])$
makes this algebra into a \dg-algebra, to be denoted
$\REnd^\bullet_{_{B}}(M):=
\bigoplus_{n\in\Z}\,\Hom'(P, P[n]).$ 

Let $\op{DGAlg}$ be the category obtained from the category of
dg-algebras and \dg-algebra morphisms by
localizing at quasi-isomorphisms.
The \dg-algebra $\REnd^\bullet_{_{B}}(M)$ viewed as an object of
 $\op{DGAlg}$  does not depend
on the choice of  projective representative $P$.
More precisely, let $\op{QIso(B)}$ denote the groupoid
that has the same objects as the category $D(B)$ and 
whose morphisms are the isomorphisms in $D(B)$. Then, one can show,
cf. \cite{Hi} for a similar result, that associating
to $M\in D(B)$ the \dg-algebra $\REnd^\bullet_{_{B}}(M)$ gives
a well-defined functor $\op{QIso(B)}\to \op{DGAlg}$.

The lift $\rdeform: \sym(T[-2])\to \RHom_{\aa\bimod}(\aa,\aa):=
\REnd_{\aa\o\aa^\opp}(\aa),$
whose existence is stated in Theorem \ref{General_formality1},
should be understood as a morphism in $\op{DGAlg}$.

For any \dg-algebra
morphism
$f: B_1\to B_2$, we 
let $f_*:
D(B_1)\to D(B_2)$ be
the push-forward
functor
$M\mapsto B_2\lo_{B_1} M,$
and $f^*: D(B_2)\to D(B_1)
$ the  pull-back functor,
 given by the change of scalars.
The functor $f^*$ is clearly $t$-exact;
it is the right adjoint of $f_*$.
These functors  are triangulated equivalences
quasi-inverse to each other, provided
the map $f$ is a \dg-algebra quasi-isomorphism.

\subsection{Homological algebra associated with a deformation.}\label{exact}
Let $\oo$ be
a formally smooth complete local algebra with maximal ideal $\fm$.
We fix   a $\k$-algebra $\aa$ and let $A$
be an infinite order formal $\oo$-deformation  of $\aa$,
as in \S\ref{inf}. Note that $A$  is a {\em flat} $\oo$-algebra.
Associated with $A$ and $\aa$, we have the corresponding
ideals $\I_{A}\sset A\o A$ and $\I_\aa\sset\aa\o\aa,$ respectively.

Set $T:=(\m/\m^2)^*$.
The projection $\oo\onto\oo/\fm^2$ induces
an isomorphism $\Tor^\oo_1(\k,\k)\iso$
$\Tor^{\oo/\fm^2}_1(\k,\k)=T.$
It follows, since  $A$ is flat over $\oo$, that the exact sequence
in \eqref{I_A} as well as all other constructions of \S\ref{setting}
are still valid in the present setting 
of  formally smooth complete local algebras $\oo$.
In particular, we have
the canonical morphism
 $u: \k\ten_{_\oo}\I_{\AA}\ten_{_\oo}\k\to
\aa\o\aa,$ 
cf. \eqref{I_A}, and the  object
$\cone(u)\in 
 D(\aa\otimes\aa^{\operatorname{op}})$.
From  \eqref{I_A} we deduce
$H^0(\cone(u))=\aa$ and
$H^{-1}(\cone(u))=\aa\o T^*.$
So, one may view  
\eqref{I_A} as an exact triangle 
\beq{du}\Delta_u:\
 \aa\o T^*[1]=H^{-1}(\cone(u))[1]\too\cone(u)\too H^0(\cone(u))=\aa,
\eeq
with  boundary map $\pa_u : \aa\to\aa\o T^*[2]$.
In this language,
the bijection of Theorem \ref{Ger_def} 
 assigns
to a deformation $(A,\psi)$ the class
 \beq{pa}
\pa_u \in
\Hom_{D(\aa\o\aa^\opp)}\bigl(\aa,\aa\o T^*[2]\bigr)
=\Ext_{\aa\o\aa^\opp}^2(\aa,\aa)\o T^*.
\eeq

There is also a different interpretation of the triangle $\Delta_u$.
Specifically, apply
 derived tensor product functor
$D(\aa\o A^\opp)\times
 D(A\o\aa^\opp)\map D(\aa\o\aa^\opp)$ to
 $\aa,$ viewed as an object of 
either $D(\aa\o A^\opp)$ or
$D(A\o\aa^\opp).$ 
This way, we get an object $\aa\lo_{_{\AA}}\aa\in  D(\aa\o\aa^\opp).$

\begin{prop}\label{C=p}
\vi The  object $\aa\lo_{_{\AA}}\aa\in  D(\aa\o\aa^\opp)$
is concentrated in {\em non-positive} degrees, and
 one has a natural quasi-isomorphism
$\dis\phi:\ta\qisto\cone(u),$ such that the
following diagram commutes
$$
\xymatrix{
{\underset{}{\ta}}\ar[d]_<>(.5){\op{qis}}^<>(.5){\phi}
\ar@{->>}[rr]^<>(.5){\op{proj}}&&
H^0(\underset{}{\aa\lo_{_{\AA}}\aa})\ar@{=}[r]&\aa\ar@{=}[d]^<>(.5){\id_\aa}\\
\cone(u)\ar@{->>}[rr]^<>(.5){\op{proj}}&&
H^0(\cone(u))\ar@{=}[r]&\aa
}
$$

\vii Thus, associated with  a deformation $(A,\psi)$
we have a canonical exact triangle
$$
\Delta_{A,\psi}:\ \aa\o T^*[1]\map\ta\map\aa,
$$
cf. \eqref{du}, with boundary map $\pa_{A,\psi},$ and
the bijection of Theorem \ref{Ger_def} reads
 \beq{ddel}
(A,\psi)\mto {\mathsf{deform}}(A,\psi)=
\pa_{A,\psi}\in 
\Ext_{\aa\o\aa^\opp}^2(\aa,\aa\o T^*).
\eeq
\end{prop}
\proof Since $A$ is flat over $\oo$,
on the category of left $A$-modules one has an isomorphism
of functors $\aa\lo_{_A}(-)=\k\lo_\oo(-)$.
Now, use \eqref{classI} to  replace the second tensor
factor $\aa$ in $\aa\lo_{_{\AA}}\aa$ by
 $\cone\big[(\I_{A}\otimes_{_{\oo}}\k)[1]
\to A\otimes \aa\big],$ a quasi-isomorphic object.
We find
$$
\aa\lo_{_{A}}\aa=\k\lo_\oo\aa=
\k\lo_\oo\cone\big[(\I_{A}\otimes_{_{\oo}}\k)[1]
\to A\otimes \aa\big]=\cone\big[\k\lo_{_\oo}(\I_{A}\otimes_{_\oo}\k)[1]
\to\aa\otimes \aa\big].
$$

The object  $\k\lo_{_\oo}(\I_{A}\otimes_{_\oo}\k)$
is concentrated in non-positive degrees,
and we clearly have
$$\bigl(\k\lo_{_\oo}(\I_{A}\otimes_{_\oo}\k)\bigr)^{\taug 0}=
H^0(\k\lo_{_\oo}(\I_{A}\otimes_{_\oo}\k))=
\k\o_{_\oo}\I_{A}\otimes_{_\oo}\k.
$$
Thus, we conclude that the object $\ta$ is quasi-isomorphic to
\begin{align*}
\cone\big[\k\lo_{_\oo}(\I_{A}\otimes_{_\oo}\k)[1]
\to\aa\otimes \aa\big]^{\taug -1}
&=\cone\big[(\k\lo_{_\oo}(\I_{A}\otimes_{_\oo}\k))^{\taug 0}[1]
\to\aa\otimes \aa\big]\\
&=\cone\big[(\k\o_{_\oo}\I_{A}\otimes_{_\oo}\k)[1]
\to\aa\otimes \aa\big]=\cone(u).\quad\Box
\end{align*}

%kkk
\subsection{Koszul duality.}\label{kos_duality} Fix a finite dimensional
vector space $T$ and
let  $\La=\wedge^\bullet(T^*[1])$ 
be the exterior algebra of the dual vector space $T^*$, 
placed in degree $-1$. For each
$n=0,-1,-2,\ldots,$
we have a graded ideal $\La_{<n}\sset\La.$
One has a canonical
 extension
 of graded $\La$-modules
\beq{klak}
\Delta_\wedge:\
0\map T^*[1]
\map\La/\La_{<-1}\stackrel{\epsilon_\wedge}\map \k_\wedge\map0,
\eeq
where we set $\k_{_\wedge}:=\La/\La_{<0}.$
We will often view $\La$ as a \dg-algebra 
concentrated in nonpositive degrees,   with zero
differential.

Recall  that the standard Koszul resolution of $\k_{_\wedge}$
provides  an explicit \dg-algebra model
for  $\REnd^\bullet_{_{\La}}(\k_{_\wedge})$ together with
an imbedding
of the graded Symmetric algebra $\sym(T[-2])$
as a subalgebra of cocycles
in that \dg-algebra model.
Furthermore,   $\La=\wedge^\bullet(T^*[1])$ is a Koszul algebra,
cf. \cite{BGG}, \cite{GKM}, so
this 
imbedding induces a graded algebra isomorphism on cohomology:
\beq{v}
{\mathsf{koszul}}:\
\sym(T[-2])\iso \Ext^\bullet_\La(\k_\wedge,\k_\wedge).
\eeq
Thus, the imbedding yields a \dg-algebra  quasi-isomorphism
\beq{kos}
{\mathfrak{Koszul}}:\
\sym(T[-2])\to \REnd^\bullet_{_{\La}}(\k_{_\wedge}),
\eeq
provided the graded algebra $\sym(T[-2])$
is viewed as a \dg-algebra with zero differential.

From \eqref{v}, we get a canonical vector space isomorphism
$$
\xymatrix{
\End_\k T=T\o T^*
\ar[rrr]^<>(0.5){{\mathsf{koszul}}\o\id_{T^*}}&&&
\Ext^2_\La(\k_\wedge,\k_\wedge)\o T^*=
\Ext^2_\La(\k_\wedge,T^*).
}
$$
It is immediate from the definition of the Koszul complex
that  the above  isomorphism
sends the element $\id_T\in \End_\k T$ to 
$\pa_\wedge\in
\Ext^2_\La(\k_\wedge,T^*),$
the class of the boundary map
$\k_\wedge\to T^*[2]$ in the canonical exact
triangle $\Delta_\wedge$ given by \eqref{klak}.

\subsection{A \dg-algebra.} 
Let $\oo=\k[[T]]$ be  the algebra of
 formal power series,  with maximal ideal $\m\sset\oo$ such that 
$T=(\fm/\fm^2)^*.$
There is a standard super-commutative
 \dg-algebra $\r$ over $\oo$
concentrated in non-positive degrees and such that
\begin{equation}\label{Ka_conditions}
\begin{array}{lll}
\vi& \r=\bigoplus_{i\leq 0}\,\r^i & \text{and}\quad\r^0=\oo,\\
\vii& H^0(\r)=\k& \text{and}\quad H^i(\r)=0\;,\;\forall i\leq -1,\\
\viii& \r\enspace\,\text{is a free graded}&\oo\text{-module.}
\end{array}
\end{equation}

To construct $\r$,
 for each $i=0,1,\ldots,n=\dim T,$
we let $\r^{-i}=\k[[T]]\o \wedge^i T^*$ be the $\oo$-module
of differential forms on the scheme
$\Spec\oo$. We put $\r:=\bigoplus_{-n\leq i\leq 0}\r^i$.
Further, write $\xi$ for the {\em Euler vector field}
on $T$. Contraction with $\xi$ gives
a differential $d: \r^{-i}\to \r^{-i+1}$,
and it is well-known that the resulting \dg-algebra
is acyclic in negative degrees, i.e.,  property 
\eqref{Ka_conditions}(ii) holds true.  
Properties \eqref{Ka_conditions}(i) and (iii) are
clear.\qed

Until the end of this section, we will use the convention that
each time a copy of 
the vector space $T^*$ occurs in a formula,
this copy has grade degree $-1$.
We form the \dg-algebra
$\dis\r\otimes_{_{\oo}}\r\simeq \k[[T]]\o \wedge^\bullet(T^*\oplus
T^*)$.
 Let $\r_\Delta\sset\r\otimes_{_{\oo}}\r$ be
the $\oo$-subalgebra generated by  the diagonal copy
$T^*\sset T^*\oplus T^*=\wedge^1(T^*\oplus T^*)$.

\begin{lemma}\label{R}
There is a   \dg-algebra imbedding
$\imath:\La\into {\r\otimes_{_{\oo}}\r}$ such that

\vi Multiplication in $\r\otimes_{_{\oo}}\r$ induces a \dg-algebra isomorphism
$\r_\Delta\o\imath(\La)\iso{\r\otimes_{_{\oo}}\r}.$

\vii The  kernel of multiplication map
$\bm_\r: {\r\otimes_{_{\oo}}\r}\to\r$ is the
ideal in the algebra ${\r\otimes_{_{\oo}}\r}$
generated by $\imath(\La_{<0})$.
\end{lemma}
\begin{proof}  
We have
 $\r\otimes_{\oo}\r\simeq \k[[T]]\o \wedge(T^*\oplus T^*)
\simeq \r_\Delta\o\wedge(T^*),$
where the last factor $\wedge(T^*)$ is generated by
the  anti-diagonal copy
$T^*\subset T^*\oplus T^*$.
It is clear that this anti-diagonal copy of
$T^*$ is
annihilated by the differential
in the \dg-algebra $\r\otimes_{\oo}\r$.
We deduce that the subalgebra generated
by the anti-diagonal copy of
$T^*$ is isomorphic to $\La$ as a \dg-algebra.
This immediately implies properties (i)-(ii).
\end{proof}

Now, let $\oo$ be an arbitrary smooth complete
local algebra.
A pair $(\r,\eta),$ where  $\r=\bigoplus_{i\leq 0}\r^i$
is a super-commutative \dg-algebra
 concentrated in non-positive degrees and 
$\eta: \oo\to \r^0$ is
an algebra homomorphism, will be referred to
as an
 \odg-algebra.
A map $h: \r\to\r',$ between two \odg-algebras
$(\r,\eta)$ and $(\r',\eta'),$
is said to be an \odg-algebra morphism
if it is a \dg-algebra map such that
$h\ccirc \eta=\eta'.$

Let $D$ denote the standard `de Rham  dg-algebra
of the line',
that is, a free supercommutative $\oo$-algebra
with one even generator $t$ of degree $-2$ and 
one odd generator
$dt$ of degree $-1$, and equipped with the $\oo$-linear differential
sending $t$ to $dt$. For any  $z\in \k$, the assignment
$t\mapsto z, \,dt\mapsto 0$ gives an \odg-algebra morphism
$\pr_z: D\onto \oo$, where $\oo$ is viewed as a \dg-algebra
with zero differential.
A pair $h,g:\r\to\r',$ of \odg-algebra morphisms,
is said to be {\em homotopic}
\footnote{The reader is referred to \cite{BoGu} for an excellent
exposition of the homotopy theory of \dg-algebras.}
provided there exists an \odg-algebra morphism
${\mathbf{h}}: \r\to D\o_\oo \r'$ such that 
the composite
$\r\stackrel{\mathbf{h}}\too D\o_\oo \r'\stackrel{\pr_z\o\id}\too\oo\o_\oo\r'=\r'$
is equal to $h$ for $z=0$, resp.  equal to $g$ for $z=1$.
Let $\op{DGCom}(\oo)$ be
the category whose objects are
 \odg-algebras and whose morphisms
are obtained 
from homotopy classes of  \odg-algebra
morphisms by localizing at  quasi-isomorphisms.

Write $\m$ for the maximal ideal in $\oo$ and put $T:=(\m/\m^2)^*$.
So we can use all the previously introduced
notation, such as $\La=\wedge^\bullet(T^*)$ (with $T^*$ in degree $-1$).

\begin{lemma}\label{homotopy} \vi  There exists an  \odg-algebra $\r$
satisfying the three conditions in \eqref{Ka_conditions}
and such that all the statements of Lemma \ref{R} hold.

Furthermore, let $\r_s,\, s=1,2,$ be two such \odg-algebras. Then, 
for $s=1,2,$we have

\vii There  exists a third \odg-algebra $\r,$ as in \vi\!\!,
and    \odg-algebra morphisms
 $h_s: \r\qisto\r_s;$ these morphisms are unique up to
 homotopy. 

Thus, the object of $\op{DGCom}(\oo)$ arising
from any choice of  \dg-algebra $\r,$ as in \vi\!\!,
is uniquely determined up to a canonical (quasi)-isomorphism.

\viii Let $\imath:\La\into {\r\otimes_{_{\oo}}\r}$ and
$\imath_s:\La\into {\r_s\otimes_{_{\oo}}\r_s},$
be the corresponding maps of 
Lemma \ref{R}\vi\!\!. 
Then, the \dg-algebra morphism
$(h_s\o h_s)\ccirc\imath$ is homotopic to $\imath_s.$
\end{lemma}

\begin{proof} Any choice of representatives in $\m$ of some
basis of the vector space $T^*=\m/\m^2$
provides a topological algebra isomorphism
$\oo\cong\k[[T]]$. This proves (i).

To prove (ii), choose an identification
$\oo\cong\k[[T]]$ and let $\r:=\k[[T]]\o \wedge^\bullet T^*$ be the 
corresponding standard \dg-algebra constructed earlier. 
Since $\r^1$ is a free $\oo$-module,
we may
find  $\oo$-module maps $h^1_s,\,s=1,2,$ that make the following
diagrams commute
$$
\xymatrix{
\r^1\ar[d]^<>(0.5){h^1_s}\ar[r]^<>(0.5){d}&\r^0=\oo\ar@{=}[d]^<>(0.5){\id}
\ar[r]^<>(0.5){\epsilon_{_\r}}&\k\ar@{=}[d]^<>(0.5){\id}\\
\r^1_s\ar[r]^<>(0.5){d}&\r^0_s=\oo\ar[r]^<>(0.5){\epsilon_{_{\r_s}}}&\k.
}
$$

Further,  $\r$ is free as a 
super-commutative $\oo$-algebra. Hence, the $\oo$-module
map $h^1_s: \r^1\to\r^1_s$ can be uniquely extended,
by multiplicativity, to a graded algebra
map $h_s:\r\to\r_s$. The latter map automatically
commutes with the differentials and, moreover, induces 
isomorphisms on cohomology, since each
algebra has no cohomology in degrees $\neq 0$.
This proves the existence of quasi-isomorphisms.
The remaining statements involving homotopies are proved 
similarly.
\end{proof}

\section{Proofs} 
\subsection{}
Fix  an associative
algebra $\aa$,  a complete smooth local $\k$-algebra $\oo$
with maximal ideal $\m$,
and an \odg-algebra $\r,$ as in Lemma \ref{homotopy}(i).
Let $(A,\psi)$ be an infinite order
formal $\oo$-deformation  of~$\aa.$

The   differential and the grading on $\r$
 make the
tensor product $\ra:=\r\otimes_{_{\oo}}{A}$  a \dg-algebra 
which is concentrated in nonpositive degrees and is such
that the subalgebra $A=A\o 1\sset\ra$ is placed in degree zero.
Since ${A}$ is flat over $\oo $,
one has $H^\bullet(\r\otimes_{_{\oo}}{A})=$
$H^\bullet(\r)\otimes_{_{\oo}}{A}=\k\otimes_{_{\oo}}{A}.$
Thus, we have a natural projection
\beq{p}
p:\
\ra\onto H^0(\ra)=\k\otimes_{_{\oo}}{A}=A/\m A\stackrel{\psi}\iso\aa
\eeq
The map $p$ is a  \dg-algebra quasi-isomorphism
that makes  $\ra$ a \dg-{\em algebra resolution} of $\aa$.
In particular, we have  mutually
quasi-inverse   equivalences
$\dis p^*,p_*: D(\ra\o\ra^\opp)\leftrightarrows
D(\aa\o\aa^\opp).$
The first surjection in \eqref{p}
may be described alternatively as the map $\epsilon_\r\o\id_A:
\r\otimes_{_{\oo}}{A}\to \k\otimes_{_{\oo}}{A},$
where  $\epsilon_\r: \r\onto \r/\r_{<0}\cong
\r^0=\oo\onto\oo/\m=\k$ is the  natural  augmentation that makes
$\r$ a \dg-algebra resolution of $\k$.

The \dg-algebra $\ror$ is also concentrated in nonpositive degrees.
We have 
\beq{rork}
(\r\otimes_{_{\oo}}\r)\lo_\La\k_{_\wedge}\iso
H^0\bigl((\r\otimes_{_{\oo}}\r)\lo_\La\k_{_\wedge}\bigr)
=
(\r\otimes_{_{\oo}}\r)\o_\La\k_{_\wedge}
\stackrel{_{\bm_\r}}\iso {\r},
\eeq
Here, the last isomorphism is
obtained by applying the functor $(-)\o_\La\k_{_\wedge}$
to  the  isomorphism $\ror\cong
\r_\Delta\o\La,$ provided by Lemmas  \ref{R}(i)-\ref{homotopy}(i).
Similarly,  applying the functor $(-)\o_\oo A$, we obtain
\beq{ror}
(\ror)\o_\oo A\iso(\r_\Delta\o\La)\o_\oo A\iso(\r_\Delta\o_\oo A)\o\La\iso
\ra\o\La.
\eeq
Let $\xi$ denote the composite isomorphism in  \eqref{ror}.
The isomorphisms \eqref{ror}, resp. \eqref{rork},
 are incorporated in the top, resp. bottom,
row of the  following natural commutative diagram of \dg-algebra maps
\beq{1}
\xymatrix{
\ra\o_A\ra\ar@{=}[r]\ar[d]_<>(.5){\bm_{_\ra}}&
(\r\o_\oo A)\o_A(\r\o_\oo A)
\ar@{=}[r]^<>(.5){A\o_A A=A}\ar[d]_<>(.5){\bm_\r\o\bm_A}&
(\ror)\o_\oo A
\ar@{=}[r]^<>(.5){\xi}\ar[d]_<>(.5){\id_{\ror}\o\epsilon_\wedge}&
\ra\o\La\ar[d]_<>(.5){\id_{\ra}\o\epsilon_\wedge}\\
\ra\ar@{=}[r]&\r\o_\oo A\ar@{=}[r]_<>(.5){\eqref{rork}\o\id_A}&
{\bigl((\ror)\o_\La\k_\wedge\bigr)\o_\oo A\,}
\ar@{=}[r]_<>(.5){\xi\o_\La\k_\wedge}&\ra\o\k_\wedge.
}
\eeq
Using the isomorphisms in the top row,
 we may (and will) further identify $(\ror)\o_\oo A$
with $\ra\o_A\ra$.
In particular,
 any $(\ror)\o_\oo A$-module may be viewed
as an $\ra$-bimodule, and we may also
view $\ra$ and $\La$ as  \dg-subalgebras in the
 \dg-algebra $\ra\o_A\ra$.

The algebra $\r$ being graded-commutative,
any graded  left $\r$-module may be
also viewed as
a right $\r$-module. Thus, any graded  left
 $\ror$-module may be viewed as
a graded  $\r$-bimodule.
A key role in our proof of Theorem
\ref{General_formality1} will be played
by  the following push-forward functor 
\begin{equation}\label{Theta}
\Th :
D(\La)\to D(\ra\o\ra^\opp),\enspace M\mapsto\bigl((\ror)\o_\oo
A\bigr)\o_\La M=
(\ra\o_A\ra)\o_\La M.
\end{equation}

By \eqref{ror}, for any \dg-module $M$ over $\La$, we get
$$
\Th (M)=\bigl((\ror)\o_\oo A\bigr)\o_\La M\simeq (\ra\o\La)\o_\La
M=\ra\o M.
$$
Although this isomorphism does {\em not} exhibit  the
$\ra$-bimodule structure on the object on the right,
it does imply that  formula \eqref{Theta} gives a well-defined
 triangulated functor. Moreover, this functor  is
 $t$-{\em exact}; indeed, since
$\ra$ is quasi-isomorphic
to $\aa$, we find 
\beq{dgv}
H^\bullet(\Th (M))\cong H^\bullet(\ra)\o
H^\bullet(M)\cong \aa\o H^\bullet(M).
\eeq

\begin{prop}
\label{2} 
For any infinite order deformation $(A,\psi)$, in $D(\aa\o\aa^\opp)$,
there is
a natural isomorphism $f_{A,\psi}: p_*\ccirc\Th (\k_{_\wedge}) \iso\aa$
  that makes
the following diagram commute
\begin{equation}\label{dia}
\xymatrix{
\sym(T[-2])\ar[rr]_<>(.5){\sim}^<>(.5){{\mathsf{koszul}}}
\ar[d]^{\sf{deform}}&&
\Ext_{_{\La}}^\bullet(\k_\wedge,\k_\wedge)
\ar[r]^<>(.5){\Th  }&
\Ext_{_{\rop}}^\bullet\bigl(\Th (\k_{_\wedge}),\Th (\k_{_\wedge})\bigr)
\ar[d]_<>(.5){p_*}^<>(.5){\op{qis}}\\
\Ext_{_{\aa\otimes\aa^\opp}}^\bullet(\aa,\aa)&&&
\Ext_{_{\aa\otimes\aa^\opp}}^\bullet\bigl(p_*\ccirc\Th (\k_{_\wedge}),
p_*\ccirc\Th (\k_{_\wedge})\bigr).
\ar[lll]_<>(.5){f_{A,\psi}}^<>(.5){\sim}
}
\end{equation}
\end{prop}
 
\begin{proof} First of all, using the definition of $\Th$
and the isomorphisms in the bottom row of   diagram \eqref{1},
we find
\beq{Thk}
\Th (\k_\wedge)=
\bigl((\ror)\o_\oo A\bigr)\o_\La \k_\wedge
=\ra.
\eeq 
Write $g: 
\Th (\k_\wedge)\iso \ra$ for the composite, and
${\mathsf{can}}_p: \ra\qisto p^*\aa$ for the map $p$ viewed
as
a morphism in $D(\ra\o\ra^\opp).$ We define a morphism
$f_{A,\psi}$ to be the following composite
$$f_{A,\psi}:\
\xymatrix{
p_*\ccirc\Th(\k_\wedge)
\ar[rr]^<>(0.5){p_*(g)}_<>(0.5){\sim}&&
p_*(\ra)\ar[rr]^<>(0.5){p_*({\mathsf{can}}_p)}_<>(0.5){\sim}&&
p_*(p^*\aa)
\ar[rr]^<>(0.5){\text{adjunction}}_<>(0.5){\sim}&&
\aa
}
$$

We claim that, with this definition of  $f_{A,\psi}$,
 diagram \eqref{dia} commutes.
To see this, observe first that all the maps in the diagram
 are clearly algebra
homomorphisms. Hence, it suffices to verify commutativity
of \eqref{dia}
on   the generators of the algebra 
$\sym(T[-2])$, that is, we must prove 
that for all $t\in T[-2]$ one has
${\sf{deform}}(t)=f_{A,\psi}\ccirc p_*\ccirc\Th  
\ccirc{{{\mathsf{koszul}}}}(t).$ 

It will be convenient to work `universally' over $T$;
that is, we treat the map ${\sf{deform}}: T[-2]\to \Ext^2_{\aa\o\aa^\opp}(\aa,\aa)$
as an element ${\sf{deform}}(A,\psi)\in \Ext^2_{\aa\o\aa^\opp}(\aa,\aa\o T^*)$,
see \eqref{ddel}. We also have the element
$\id_T\in \End_\k T=\sym^1(T[-2])\o T^*$.
Now, tensoring with
$T^*$, we rewrite the equation
that we must prove as 
$\,\dis{\sf{deform}}(A,\psi)=
f_{A,\psi}\ccirc p_*\ccirc\Th  
\ccirc{\mathsf{koszul}}(\id_T).$  Both sides here belong to
$ \Ext^2_{\aa\o\aa^\opp}(\aa,\aa\o T^*).$
Thus, applying further  $p^*(-)$ to each side and using  adjunctions, 
we see that proving the
Proposition amounts to 
showing that
\begin{equation}\label{prove1}
p^*({\sf{deform}}(A,\psi))\;=\;
\Th \ccirc{\mathsf{koszul}}(\id_T)
\quad\text{holds in}\quad
\Ext^2_{\ra\o\ra^\opp}(\ra,\ra\o T^*).
\end{equation}

We compute the LHS of this equation using Proposition  \ref{C=p}(ii),
which says  ${\sf{deform}}(A,\psi)=\pa_{A,\psi}.$
Therefore, $p^*({\sf{deform}}(A,\psi))=
p^*(\pa_{A,\psi}),$ is the boundary map for
 $p^*(\Delta_{A,\psi})$,
 the pull-back via the equivalence  $p^*:{D(\aa\o\aa^\opp)}\iso
D(\ra\o\ra^\opp)$ of the canonical triangle 
$\Delta_{A,\psi}$ 
that appears in part (ii) of Proposition \ref{C=p}.
Now, using the quasi-isomorphism
$p^*\aa\cong\ra,$ we can write the triangle  $p^*(\Delta_{A,\psi})$
as follows
$$
p^*(\Delta_{A,\psi}):\
\ra\o T^*[1]\tooo p^*\bigl((\aa\lo_A\aa)^{\taug -1}\bigr)
\stackrel{p^*(\bm_\aa)}\tooo\ra.
$$

To describe the middle term in the  last triangle we recall
that $\ra$ is {\em free} over the subalgebra $A\sset \ra$.
It follows that in $D(\ra\o\ra^\opp)$ one has
$p^*(\aa\lo_A\aa)\cong\ra\lo_A\ra\cong\ra\o_A\ra.$
Hence, we deduce
$p^*\bigl((\aa\lo_A\aa)^{\taug -1}\bigr)=(\ra\o_A\ra)^{\taug -1},$
since  the pull-back functor $p^*$ is  always $t$-exact.
 Thus, we see that our exact
triangle takes the following final form
\beq{LHS}
p^*(\Delta_{A,\psi}):\
\ra\o T^*[1]\tooo (\ra\o_A\ra)^{\taug -1}
\stackrel{\bm_\ra}\tooo\ra.
\eeq

Next, we analyze the RHS of equation \eqref{prove1}.
By \S\ref{kos_duality},
 we have
${\mathsf{koszul}}(\id_T)=\pa_\wedge.$ Hence, the class
$\dis\Th  
\ccirc{\mathsf{koszul}}(\id_T)=\Th(\pa_\wedge),$ in 
$\Hom_{D(\ra\o\ra^\opp)}(\ra,\ra\o T^*[2]),$ is represented by the
boundary map for the exact triangle $\Th(\Delta_\wedge),$
see \eqref{klak}. The  latter reads
\beq{klak2}
\Th(\Delta_\wedge):\
 \Th(T^*[1])
\tooo\Th(\La/\La_{<-1})
\stackrel{\Th(\epsilon_\wedge)}\tooo \Th(\k_\wedge).
\eeq
Here, $\Th(\k_\wedge)=\ra$, by \eqref{Thk}, hence $\Th(T^*)=\ra\o T^*.$
Further, by definition we have
$$\Th(\La/\La_{<-1})=(\ra\o_A\ra)\o_\La (\La/\La_{<-1})=
(\ra\o_A\ra)/(\ra\o_A\ra)\cd\La_{<-1}.$$
Now, $\ra\o_A\ra\cong\ra\o\La$, by  \eqref{ror}, and we deduce
\beq{iso1}
(\ra\o_A\ra)\cd\La_{<-1}\cong\ra\o\La_{<-1}.
\eeq
We see that
the morphism $\Th(\epsilon_\wedge)$
 in \eqref{klak2} may be viewed as a map
induced 
by the leftmost vertical arrow in diagram  \eqref{1}.
Thus, the triangle in \eqref{klak2} takes the following form
(note that Lemma \ref{R}(iii) insures that $\bm_\ra$ maps
$(\ra\o_A\ra)\cd\La_{<-1}$ to zero)
\beq{RHS}
\Th(\Delta_\wedge):\
\ra\o T^*[1]\too(\ra\o_A\ra)/(\ra\o_A\ra)\cd\La_{<-1}
\stackrel{\bm_\ra}\too \ra
\eeq

To compare the LHS with the RHS of \eqref{prove1},
one has to compare \eqref{LHS} with \eqref{RHS}. We see that in order
to prove \eqref{prove1} it suffices to show 
\beq{LR}
(\ra\o_A\ra)^{\taug -1}=
(\ra\o_A\ra)/(\ra\o_A\ra)\cd\La_{<-1}
\quad\text{in}\enspace
D(\ra\o\ra^\opp).
\eeq

To this end, we
write an exact triangle
\beq{vect2}
(\ra\o_A\ra)\cd\La_{<-1}\to\ra\o_A\ra\to
(\ra\o_A\ra)/(\ra\o_A\ra)\cd\La_{<-1}.
\eeq
Here, the \dg-vector space  on the right  is isomorphic
to $\Th(\La/\La_{<-1})$, hence, has no cohomology in degrees ${<-1},$
by \eqref{dgv}. Similarly, the \dg-vector space
on the left is  isomorphic
to $\ra\o\La_{<-1}$, see \eqref{iso1}, hence, has no cohomology in
degrees  ${\geq-1}.$ Therefore, the triangle in \eqref{vect2} must be isomorphic
to the canonical exact triangle
$\dis(\ra\o_A\ra)^{\taul -1}\to \ra\o_A\ra\to (\ra\o_A\ra)^{\taug -1}$.
This proves \eqref{LR}.
\end{proof}

\subsection{Proposition \ref{2} implies  Theorem
\ref{General_formality1}.}
To see this, we observe that the functors $\Th $ and $p^*$
provide  us not only with  maps between the Ext-groups which occur
in diagram \eqref{dia}, but also with  lifts
of those maps to   morphisms in $\op{DGAlg}$:
\beq{PHI}
\xymatrix{
\REnd_{_{D(\La)}}^\bullet(\k_\wedge)
\ar[rr]^<>(0.5){p_*\ccirc\Th  }&&
\REnd_{_{D(\aa\otimes\aa^\opp)}}^\bullet\bigl((p_*\ccirc\Th (\k_{_\wedge})\bigr)
\ar[r]^<>(0.5){f_{A,\psi}}_<>(0.5){\sim}&
\REnd_{_{D(\aa\otimes\aa^\opp)}}^\bullet(\aa).
}
\eeq

Let ${\mathfrak{Deform}}:=p_*\ccirc\Th  \ccirc{\mathfrak{Koszul}}$
be the composite of the 
 \dg-algebra morphism $\sym(T[-2])\to\REnd_{_{\La}}^\bullet(\k_\wedge)$,
see \eqref{kos}, followed by the morphisms
in \eqref{PHI}. Thus, in   $\op{DGAlg}$, we get a morphism 
 ${\mathfrak{Deform}}:
\sym(T[-2])\to\REnd_{_{D(\aa\otimes\aa^\opp)}}^\bullet(\aa).$
Further,
 the induced map of cohomology
$$H^\bullet({\mathfrak{Deform}})=H^\bullet(p_*\ccirc\Th  \ccirc{\mathfrak{Koszul}})=
H^\bullet(p_*\ccirc\Th  )\ccirc{\mathsf{koszul}}:\
\,
 \sym(T[-2])\to\Ext_{_{\aa\otimes\aa^\opp}}^\bullet(\aa)$$
is equal to the map $\mathsf{deform}$,
by Proposition \ref{2}.
Thus, the  morphism
 ${\mathfrak{Deform}}$
yields a  morphism  in   $\op{DGAlg}$ as required by Theorem
\ref{General_formality1}.

Our construction of the functor $\Th$, hence of the
morphism ${\mathfrak{Deform}},$ was based on the choice of
an \odg-algebra $\r$. To show independence of such a choice,
let $\r_s,\,s=1,2,$ be two  \odg-algebras,
as in Lemma \ref{homotopy}, and write $\ra_s:=\r_s\o_{_\oo}\aa.$
Part (ii) of the Lemma
implies that there exists a canonical 
isomorphism ${h}: \ra_1\qisto\ra_2,$ in   $\op{DGAlg}$,
which is compatible with the augmentations
$p_s: \ra_s\qisto\aa$. The  isomorphism ${h}$ induces
a triangulated equivalence 
${h}_*: D(\ra_1\o\ra_1^\opp)
\iso D(\ra_2\o\ra_2^\opp),$ and also
a canonical isomorphism ${h}_{\End},$
 in   $\op{DGAlg}$, between
the two \dg-algebra models
for
the object $\REnd_{_{D(\aa\otimes\aa^\opp)}}^\bullet(\aa)\in
\op{DGAlg}$,
constructed using $\r_1$ and $\r_2,$ respectively.

Now, let $\Theta_s: D(\La)\to D(\ra_s\o\ra_s^\opp),\, s=1,2,$
be the corresponding two functors defined as in \eqref{Theta}.
Lemma \ref{homotopy}(iii) yields a  canonical
 isomorphism of functors $\Phi: \Theta_2\iso {h}_*\ccirc\Theta_1.$
This way,  in  $\op{DGAlg}$, we obtain the following 
isomorphisms
$$
\xymatrix{
\REnd\bigl((p_1)_*\ccirc\Th_1
(\k_{_\wedge})\bigr)
\ar[r]^<>(0.5){{h}_*}_<>(0.5){\sim}&
\REnd\bigl({h}_*\ccirc(p_1)_*\ccirc\Th_1(\k_{_\wedge})\bigr)
\ar[r]^<>(0.5){\Phi}_<>(0.5){\sim}&
\REnd\bigl((p_2)_*\ccirc\Th_2(\k_{_\wedge})\bigr)
}
$$
where we have used shorthand notation $\REnd=\REnd_{_{D(\aa\otimes\aa^\opp)}}^\bullet.$
Let $\Upsilon$ denote the composite isomorphism,
and write $f_{A,\psi,s}$  for the  isomorphism
of Proposition  \ref{2} corresponding to
the \dg-algebra  $\r_s,\,s=1,2$.
It is straightforward to check that our construction
insures  commutativity of the following diagram  in  $\op{DGAlg}$
$$
\xymatrix{
\REnd_{_{D(\La)}}^\bullet(\k_\wedge)\ar@{=}[d]^<>(0.5){\id}
\ar[rr]^<>(0.5){(p_1)_*\ccirc\Th_1  }&&
\REnd_{_{D(\aa\otimes\aa^\opp)}}^\bullet\bigl((p_1)_*\ccirc\Th_1 (\k_{_\wedge})\bigr)
\ar[rr]^<>(0.5){f_{A,\psi,1}}\ar[d]^<>(0.5){\Upsilon}&&
\REnd_{_{D(\aa\otimes\aa^\opp)}}^\bullet(\aa)\ar[d]^<>(0.5){{h}_{\End}}\\
\REnd_{_{D(\La)}}^\bullet(\k_\wedge)
\ar[rr]^<>(0.5){(p_2)_*\ccirc\Th_2  }&&
\REnd_{_{D(\aa\otimes\aa^\opp)}}^\bullet\bigl((p_2)_*\ccirc\Th_2(\k_{_\wedge})\bigr)
\ar[rr]^<>(0.5){f_{A,\psi,2}}&&
\REnd_{_{D(\aa\otimes\aa^\opp)}}^\bullet(\aa).
}
$$

It follows from commutativity of
the  diagram that, for the resulting two
morphisms
 ${\mathfrak{Deform}}_s:
\sym(T[-2])\to\REnd_s,\, s=1,2,$ in   $\op{DGAlg}$,
we have ${h}_{\REnd}\ccirc{\mathfrak{Deform}}_1=
{\mathfrak{Deform}}_2$. Thus, the  morphisms
${\mathfrak{Deform}}$   arising from various choices of
$\r$ all give the same    morphism in $\op{DGAlg}$.  Therefore this morphism
is
{\em canonical}. \qed

\footnotesize{
\medskip

\footnotesize{
{\bf R.B.}:  Department of Mathematics,  MIT,
77 Mass. Ave, Cambridge, MA 02139, USA;\\
\hphantom{x}\enspace\quad\quad\enspace {\tt bezrukav@math.mit.edu}}
\smallskip

\footnotesize{
{\bf V.G.}: Department of Mathematics, University of Chicago, 
Chicago IL
60637, USA;\\ 
\hphantom{x}\ab\, {\tt ginzburg@math.uchicago.edu}}

\end{document}